\documentstyle[leqno,epsf]{book}
\newcounter{supersection}[section]
\newtheorem{th}[supersection]{Theorem}
\newtheorem{lm}[supersection]{Lemma}
\newtheorem{re}[supersection]{Remark}

\newtheorem{ex}[supersection]{Example}

\def\bibname{\textbf{REFERENCES}}
\def\thebibliography#1{\paragraph*{\uppercase{\bibname}}\list
{\arabic{enumi}.}{\settowidth\labelwidth{[#1]}\leftmargin\labelwidth
\advance\leftmargin\labelsep\usecounter{enumi}}
\def\newblock{\hskip .11em plus .33em minus .07em}
\sloppy\clubpenalty4000\widowpenalty4000
\sfcode`\.=1000\relax}

\def\arc#1#2{\mathop{#1#2}\limits^{\vspace*{-0.3ex}\frown}}
\def\stop{\mbox{\footnotesize {\vrule width 6pt height 6pt}}}

\def\adraw#1#2#3{\removelastskip
\vskip 0cm plus #1cm
\goodbreak
\vskip 0cm plus -#1cm \vskip 0.3cm
\vspace*{12.0 mm}  
\hspace*{6.0 mm}   
\dimen0=#2 true bp \ifdim\dimen0>\hsize\dimen0=\hsize\fi
\epsfxsize=\dimen0 \noindent {\epsfbox{#3.eps}} \count255=\epsfurx
\advance\count255 by-\epsfllx \message{[#3.eps bp width =
\number\count255]}
\goodbreak}

\def\bdraw#1#2#3{\removelastskip
\vskip 0cm plus #1cm
\goodbreak
\vskip 0cm plus -#1cm \vskip 0.3cm
\vspace*{-50.0 mm} 
\hspace*{+72.5 mm} 
\dimen0=#2 true bp \ifdim\dimen0>\hsize\dimen0=\hsize\fi
\epsfxsize=\dimen0 \noindent {\epsfbox{#3.eps}} \count255=\epsfurx
\advance\count255 by-\epsfllx \message{[#3.eps bp width =
\number\count255]}
\goodbreak}

\def\cdraw#1#2#3{\removelastskip
\vskip 0cm plus #1cm
\goodbreak
\vskip 0cm plus -#1cm \vskip 0.3cm
\vspace*{5.0 mm} 
\hspace*{0.0 mm} 
\dimen0=#2 true bp \ifdim\dimen0>\hsize\dimen0=\hsize\fi
\epsfxsize=\dimen0 \noindent {\epsfbox{#3.eps}} \count255=\epsfurx
\advance\count255 by-\epsfllx \message{[#3.eps bp width =
\number\count255]}
\goodbreak}

\begin{document}

\thispagestyle{plain}

\noindent $\;$ 

\noindent $\;$ 

\centerline{\Large \bf THE M\"{O}BIUS-POMPE\"{I}U METRIC PROPERTY}
\footnotetext{2000 Mathematics Subject Classification: 54E35, 51M16.}
\footnotetext{Research partially supported by the MNTRS, Serbia \&
Montenegro, Grant No. 1861.}

\vspace*{5.0 mm}

\centerline{\large Branko J. Male\v sevi\' c}

\vspace*{5.0 mm}

\begin{center}
University of Belgrade, Faculty of Electrical Engineering,    \\
Bulevar kralja Aleksandra 73, Belgrade, Serbia \& Montenegro  \\
E-mail: malesevic@kiklop.etf.bg.ac.yu
\end{center}

\vspace*{5.0 mm}

\begin{center}
\parbox{25.0cc}{\scriptsize \boldmath \bf
In the paper we consider an extension of M\"{o}bius-Pompe\"{i}u
theorem of the elementary geometry over metric spaces. We
specially take into consideration Ptolemaic metric spaces.}
\end{center}

\bigskip

\noindent
\section{\Large \bf \boldmath \hspace*{-5.0 mm} The M\"{o}bius-Pompe\"{i}u theorem and metric
spaces}

\setcounter{equation}{0}

In this paper we consider the following statement of elementary geometry \cite{Mitrinovic1},
\cite{Mitrinovic2}:

\begin{th}[{\bf M\"{o}bius, Pompe\"{i}u}]
Let $ABC$ be an equilateral triangle and $M$ any point in its plane.
Then segments $MA$, $MB$ and $MC$ are sides of a triangle.
\end{th}

\noindent
Let us consider analogous problem for the metric space $(X,d)$ with at least four points.
Let $A,B,C \in X$ be three fixed points. Then, for the point $M \in X$ we suppose that
{\em a triangle can be formed} from the distances $d_1=d(M,A)$, $d_2=d(M,B)$ and
$d_3=d(M,C)$ iff the following conjunction of inequalities is true:
\begin{equation}
\label{start}
d_1 + d_2 - d_3 \geq 0 \;\;\; \mbox{and} \;\;\;
d_2 + d_3 - d_1 \geq 0 \;\;\; \mbox{and} \;\;\;
d_3 + d_1 - d_2 \geq 0.
\end{equation}
If in conjunction (\ref{start}) at least one equality is true, then we suppose that
{\em a degenerative triangle can be formed}. If in (\ref{start}) sharp inequalities are true:
\begin{equation}
\label{start2}
d_1 + d_2 - d_3 > 0 \;\;\; \mbox{and} \;\;\;
d_2 + d_3 - d_1 > 0 \;\;\; \mbox{and} \;\;\;
d_3 + d_1 - d_2 > 0,
\end{equation}
then we suppose that {\em a non-degenerative triangle can be formed}.
In this case, for the point $M$, for which the conjunction (\ref{start2})
is true, we define that {\em point have M\"{o}bius-Pompe\"{i}u metric property}.
The main subject of this paper is to determine points $M$ which do not have
M\"{o}bius-Pompe\"{i}u metric property, i.e. these points which fulfill
the following disjunction of the inequalities:
\begin{equation}
\label{restart}
d_1 + d_2 - d_3 \leq 0 \;\;\; \mbox{or} \;\;\;
d_2 + d_3 - d_1 \leq 0 \;\;\; \mbox{or} \;\;\;
d_3 + d_1 - d_2 \leq 0.
\end{equation}
Let us notice that the point $M \!\in\! \{A,B,C\}$ do not have M\"{o}bius-Pompe\"{i}u
metric property. Thus in consideration which follows, we assume that the metric
space $(X,d)$ has at least four points.

\break

\bigskip

\bigskip

\noindent
\section{\Large \bf \boldmath \hspace*{-5.0 mm} Ptolemaic metric spaces}

\setcounter{equation}{0}

A metric space $(X,d)$ is called {\em Ptolemaic metric space}
if  Ptolemaic inequality holds:
\begin{equation}
\label{ptolemaic}
d(x_1,x_2) d(x_3,x_4) \leq d(x_2,x_4) d(x_1,x_3) + d(x_1,x_4) d(x_2,x_3)
\end{equation}
for every $x_1, x_2, x_3, x_4 \in X$ \cite{Hasto}. A normed space $(X,|\,.\,|)$
is {\em Ptolemaic normed space} if metric space $(X,d)$ is Ptolemaic with the
distance $d(x,y) = |x-y|$. Let us notice that the following lemma is true \cite{Hasto}:

\begin{lm}
A normed space is Ptolemaic iff it is an inner product space.
\end{lm}

\noindent
We give two basic examples of Ptolemaic spaces \cite{Hasto}.

\begin{ex}
\label{Example}
{\boldmath $1^{0}.$}
The space $\mbox{\bf R}^n$ with the Euclidean metric $\mbox{\textit{\texttt{d}}}(x,y)
= \mbox{\boldmath $|$}x-y\mbox{\boldmath $|$}$ is a Ptolemaic metric space.

\noindent
{\boldmath $2^{0}.$}
The space $\mbox{\bf R}^n$ with the chordal metric
on the unit Riemann sphere
$\overline{\mbox{\textit{\texttt{d}}}
\mathstrut}(x,y) = \displaystyle \frac{2
\mbox{\boldmath $|$}x-y\mbox{\boldmath $|$}}{
\sqrt{1\!+\!\mbox{\boldmath $|$}x\mbox{\boldmath $|$}^2}
\sqrt{1\!+\!\mbox{\boldmath $|$}y\mbox{\boldmath $|$}^2}}$
is a Ptolemaic metric space.
\end{ex}

\smallskip\noindent
We will illustrate following considerations  with the previous examples
of Ptolemaic metric spaces in the case of dimension $n=2$.

\smallskip

\noindent
\section{\Large \bf \boldmath \hspace*{-5.0 mm} The main results}

\setcounter{equation}{0}

Let $(X,d)$ be a metric space. Let us fix three points $A,B,C \in X$
and form distances:
\begin{equation}
a = d(B,C),\;
b = d(C,A),\;
c = d(A,B).
\end{equation}
For any point $M \in X$ let us form distances:
\begin{equation}
d_1=d(M,A),\;
d_2=d(M,B),\;
d_3=d(M,C).
\end{equation}

\smallskip

\leftline{\large \boldmath \bf Inequality $d_2 + d_3 \leq d_1$}

\bigskip

\noindent
Let us determine a set of $M$ points of metric spaces $X$ for which the following
inequality is true:
\begin{equation}
d_2 + d_3 \leq d_1.
\end{equation}
Let us form two functions:
\begin{equation}
\label{Alpha}
\alpha_1 = \alpha_1(M) = 4 d_2^2 d_3^2 - {\big (}d_1^2 - (d_2^2 + d_3^2){\big )}^2,
\end{equation}
\begin{equation}
\label{Beta}
\beta_1 = \beta_1(M) = d_2^2 + d_3^2 - d_1^2.
\end{equation}

\break

\begin{lm}
\label{LEMA_1}
For points $A$, $B$ and $C$ inequality $\alpha_1 \leq 0$ is true.
\end{lm}
{\bf Proof.} For point $A$: $d_1 = 0$ and $\alpha_1 = -(c^2-b^2)^2 \leq 0$ are true.
Similarly, the previous inequality is true for the points $B$ and $C$.~\stop

\smallskip

\begin{ex}
Let vertices $A$, $B$, $C$ of the triangle $ABC$ in the plane $\mbox{\bf R}^2$
be given by coordinates $A(a_1,b_1)$, $B(a_2,b_2)$, $C(a_3,b_3)$ and let $M(x,y)$
be any point in its plane.

\smallskip \noindent
{\boldmath $1^{0}.$}
Let us in the plane $\mbox{\bf R}^2$ use Euclidean metric $\mbox{\textit{\texttt{d}}}$.
Let us specify the form of term $\alpha_1$ and $\beta_1$ which correspond to functions
$(\ref{Alpha})$ and $(\ref{Beta})$ respectively. It is true$:$
\begin{equation}
\label{curve11} \;\;\;\;\;\; \alpha_1 = \mbox{\tt k}
(x^2\!+\!y^2)^2\!+\!(\mbox{\tt A}_1x\!+\!\mbox{\tt B}_1y)(x^2\!+\!y^2)
\!+\!\mbox{\tt C}_1x^2\!+\!\mbox{\tt D}_1xy\!+\!\mbox{\tt E}_1y^2
\!+\!\mbox{\tt F}_1x\!+\!\mbox{\tt G}_1y\!+\!\mbox{\tt H}_1,
\end{equation}
for some coefficients $\mbox{\tt k}, \mbox{\tt A}_1, \mbox{\tt B}_1, \mbox{\tt C}_1,
\mbox{\tt D}_1, \mbox{\tt E}_1,\mbox{\tt F}_1,\mbox{\tt G}_1,\mbox{\tt H}_1 \!\in\!
\mbox{\bf R}$ $(\mbox{\tt k} \!=\! 3)$. Equality $\alpha_1 \!=\! 0$ determines the algebraic
curve of the fourth order. By inequality $\alpha_1 \!<\! 0$ we determine the interior of the
previous curve. Also, it is true$:$
\begin{equation}
\label{curve12}
\beta_1
=
\mbox{\tt A}_2(x^2\!+\!y^2)\!+\!\mbox{\tt B}_2x\!+\!\mbox{\tt C}_2y\!+\!\mbox{\tt D}_2,
\end{equation}
for some coefficients $\mbox{\tt A}_2,\mbox{\tt B}_2,\mbox{\tt C}_2,\mbox{\tt D}_2
\!\in\! \mbox{\bf R}$ $(\mbox{\tt A}_2 \!=\! 1)$. If $\mbox{\tt B}_2^2\!+\!\mbox{\tt C}_2^2
\!>\! 4 \mbox{\tt D}_2$ equality $\beta_1 \!=\! 0$ is possible and determines the
circle. Then by inequality $\beta_1 \!<\! 0$ we determine the interior of the circle.

\medskip
\noindent
{\boldmath $2^{0}.$}
Let us in the plane $\mbox{\bf R}^2$ use chordal metric $\overline{\mbox{\textit{\texttt{d}}}
\mathstrut}$. Let us specify the form of the term $\overline{\alpha\mathstrut}_1$ and
$\overline{\beta\mathstrut}_1$ which correspond to functions $(\ref{Alpha})$ and
$(\ref{Beta})$ respectively. It is true$:$
\begin{equation}
\label{curve21}
\;\;\;\;\;\; \overline{\alpha\mathstrut}_1 =
\displaystyle\frac{
\overline{\mbox{\tt k}\mathstrut}
(x^2\!+\!y^2)^2
\!+\!
(\overline{\mbox{\tt A}\mathstrut}_1x
\!+\!
\overline{\mbox{\tt B}\mathstrut}_1y)(x^2\!+\!y^2)
\!+\!
\overline{\mbox{\tt C}\mathstrut}_1x^2
\!+\!
\overline{\mbox{\tt D}\mathstrut}_1xy
\!+\!
\overline{\mbox{\tt E}\mathstrut}_1y^2
\!+\!
\overline{\mbox{\tt F}\mathstrut}_1x
\!+\!
\overline{\mbox{\tt G}\mathstrut}_1y
\!+\!
\overline{\mbox{\tt H}\mathstrut}_1}{
(1+x^2+y^2)^2(1+a_1^2+b_1^2)^2(1+a_2^2+b_2^2)^2(1+a_3^2+b_3^2)^2},
\end{equation}
for some coefficients $
\overline{\mbox{\tt k}\mathstrut},
\overline{\mbox{\tt A}\mathstrut}_1,
\overline{\mbox{\tt B}\mathstrut}_1,\overline{\mbox{\tt C}\mathstrut}_1,
\overline{\mbox{\tt D}\mathstrut}_1,\overline{\mbox{\tt E}\mathstrut}_1,
\overline{\mbox{\tt F}\mathstrut}_1,\overline{\mbox{\tt G}\mathstrut}_1,
\overline{\mbox{\tt H}\mathstrut}_1 \!\in\! \mbox{\bf R}$.
If $\overline{\mbox{\tt k}\mathstrut} \neq 0$
equality $\overline{\alpha\mathstrut}_1 \!=\! 0$ determines
the algebraic curve of the fourth order.
Then
by inequality $\overline{\alpha\mathstrut}_1 \!<\! 0$
we determine the interior of the previous curve. Also, it is true$:$
\begin{equation}
\label{curve22}
\overline{\beta\mathstrut}_1
=
\displaystyle\frac{\overline{\mbox{\tt A}\mathstrut}_2(x^2+y^2)
\!+\!
\overline{\mbox{\tt B}\mathstrut}_2x
\!+\!
\overline{\mbox{\tt C}\mathstrut}_2y
\!+\!
\overline{D\mathstrut}_2}{
(1+x^2+y^2)(1+a_1^2+b_1^2)(1+a_2^2+b_2^2)(1+a_3^2+b_3^2)},
\end{equation}
for some coefficients $\overline{\mbox{\tt A}\mathstrut}_2,
\overline{\mbox{\tt B}\mathstrut}_2,\overline{\mbox{\tt C}\mathstrut}_2,
\overline{\mbox{\tt D}\mathstrut}_2 \!\in\! \mbox{\bf R}$.
If $\overline{\mbox{\tt A}\mathstrut}_2 \neq 0$ and
$\overline{\mbox{\tt B}\mathstrut}_2^2 \!+\! \overline{\mbox{\tt C}\mathstrut}_2^2
\!>\! 4 \overline{\mbox{\tt A}\mathstrut}_2 \overline{\mbox{\tt D}\mathstrut}_2$
equality $\overline{\beta\mathstrut}_1 \!=\! 0$ is possible and determines
the circle. Then by the inequality $\overline{\beta\mathstrut}_1 \!<\! 0$
we determine the interior of the circle.
\end{ex}

\smallskip \noindent
Further, let us notice that for the function $\alpha_1$:
\begin{equation}
\label{PROIZVOD}
\alpha_1
=
(d_2 + d_3-d_1)(d_3+d_1-d_2)(d_1+d_2-d_3)(d_1+d_2+d_3).
\end{equation}
According to (\ref{PROIZVOD}) equality $\alpha_1 = 0$ is equivalent with union of equalities:

\break

\noindent
\begin{equation}
\label{PrvaKomponenta}
\alpha_1^{(1)} = d_2 + d_3 - d_1 = 0,
\end{equation}
\begin{equation}
\alpha_1^{(2)} = d_3 + d_1 - d_2 = 0,
\end{equation}
\begin{equation}
\alpha_1^{(3)} = d_1 + d_2 - d_3 = 0.
\end{equation}
Subject to our further consideration is an inequality $\alpha_1^{(1)} \! \leq 0$.

\begin{lm}
\label{LEMA_2}
{\boldmath $1^{0}\!.$} For the point $B$$:$  $d_2 + d_3 \leq d_1$ iff $c\!\geq\!a$.
{\boldmath $2^{0}\!.$} For the point $C$$:$  $d_2 + d_3 \leq d_1$ iff $b\!\geq\!a$.
\end{lm}
\begin{re}
If $a > b,\,\!c$ then for points $B$ and $C$$:$
$\alpha_1 \leq 0$ and $\alpha_{1}^{(1)} > 0$.
\end{re}

\begin{lm}
\label{LEMA_3}
If for point $M$$:$ $d_2 + d_3 \leq d_1$, then we have inequalities$:$
\begin{equation}
\label{IQ_1}
d_1 + d_2 \geq d_3, \; \mbox{where equality is true for $M = B$ and $a = c$}
\end{equation}
and
\begin{equation}
\label{IQ_2}
d_3 + d_1 \geq d_2, \; \mbox{where equality is true for $M = C$ and $a = b$}.
\end{equation}
\end{lm}
{\bf Proof.} It is true
\begin{equation}
(d_1) + d_2 - d_3 \geq (d_2 + d_3) + d_2 - d_3 = 2d_2 \geq 0.
\end{equation}
Hence, the inequality (\ref{IQ_1}) follows. Thus, the equality is true only if
$M = B$ $(d_2 = 0)$ and $a = c$. Analogously, it is true
\begin{equation}
d_3 + (d_1) - d_2 \geq d_3 + (d_2 + d_3) - d_2 = 2d_3 \geq 0.
\end{equation}
Hence, the inequality  (\ref{IQ_2}) follows. Thus, the equality is true only if
$M = C$ $(d_3 = 0)$ and $a = b$.~\stop

\begin{lm}
\label{LEMA_4}
{\boldmath $1^{0}.$} If the point $M$ fulfills $d_2 + d_3 \leq d_1$ then the following
implication is true$:$
\begin{equation}
\label{IMPLIKACIJA_1}
\alpha_1 \leq 0 \Longrightarrow \beta_1 \leq 0.
\end{equation}

\noindent
{\boldmath $2^{0}.$} If the point $M$ fulfills $d_3 + d_1 \leq d_2$ or $d_1 + d_2 \leq d_3$
then the following implication is true$:$
\begin{equation}
\label{IMPLIKACIJA_2}
\alpha_1 \leq 0 \Longrightarrow \beta_1 \geq 0.
\end{equation}
\end{lm}
{\bf Proof.} The implications (\ref{IMPLIKACIJA_1}) and (\ref{IMPLIKACIJA_2})
have the same assumptions:
\begin{equation}
\begin{array}{rcl}
\alpha_{1} &\!\!=\!\!& 4d_2^2d_3^2 - (d_1^2 - d_2^2 - d_3^2)^2 \\[1.0 ex]
           &\!\!=\!\!&
{\big (}2d_2d_3-d_1^2+d_2^2+d_3^2 {\big )}
{\big (}2d_2d_3+d_1^2-d_2^2-d_3^2 {\big )}
\leq 0,
\end{array}
\end{equation}
which follow if the following conjunction is true
\begin{equation}
\label{KONJUKCIJA_1}
{\big (}2d_2d_3-d_1^2+d_2^2+d_3^2 {\big )} \leq 0
\;\;\mbox{and}\;\;
{\big (}2d_2d_3+d_1^2-d_2^2-d_3^2 {\big )} \geq 0
\end{equation}
or the conjunction
\begin{equation}
\label{KONJUKCIJA_2}
{\big (}2d_2d_3-d_1^2+d_2^2+d_3^2 {\big )} \geq  0
\;\;\mbox{and}\;\;
{\big (}2d_2d_3+d_1^2-d_2^2-d_3^2 {\big )} \leq 0.
\end{equation}

\break

\noindent
{\boldmath $1^{0}.$} Let $d_2 + d_3 \leq d_1$ be true.
For $M=B$ or $M=C$ implication (\ref{IMPLIKACIJA_1})
is directly verified. Especially for $M=B$ and
$a=c$ or for $M=C$ and $a=b$ equality $\beta_1 = 0$ is true.
Let us assume that $M \neq B,C$ and let us assume that $\alpha_1 \leq 0$
in (\ref{IMPLIKACIJA_1}) be true. On the basis of $d_2 + d_3 \leq d_1$,
according to lemma \ref{LEMA_3} it follows that $d_1 + d_2 > d_3$ and
$d_3 + d_1 > d_2$. Therefore
\begin{equation}
\label{IQ_11}
2d_2d_3 - d_1^2 + d_2^2 + d_3^2 = (d_2 + d_3)^2 - d_1^2 \leq 0
\end{equation}
and
\begin{equation}
\label{IQ_12}
2d_2d_3+d_1^2-d_2^3-d_3^2 = (d_1-d_2+d_3)(d_1+d_2-d_3) > 0.
\end{equation}
From (\ref{IQ_11}) and (\ref{IQ_12}) we can conclude that the conjunction
(\ref{KONJUKCIJA_1}) is true and conjunction (\ref{KONJUKCIJA_2}) is not true.
From the conjunction (\ref{KONJUKCIJA_1}) it follows that $d_1^2-d_2^2-d_3^2 \geq 2d_2d_3
> d_2^2 + d_3^2 - d_1^2$ and from there $d_1^2 > d_2^2 + d_3^2$, ie. $\beta_1 < 0$.

\smallskip \noindent
{\boldmath $2^{0}.$} Let $d_3 + d_1 \leq d_2$ be true. For $M=B$ or $M=C$ implication
(\ref{IMPLIKACIJA_2}) is directly verified. Especially for $M=B$ and $a=c$ or for $M=C$
and $a=b$ equality $\beta_1 = 0$ is true. Let us assume that $M \neq B,C$
and let us assume that $\alpha_1 \leq 0$ in (\ref{IMPLIKACIJA_2}) be true.
On the basis of $d_3 + d_1 \leq d_2$, according to the lemma analogous
to lemma \ref{LEMA_3}, it follows $d_2 + d_3 > d_1$ and
$d_1 + d_2 > d_3$. Therefore
\begin{equation}
\label{IQ_21}
2d_2d_3 - d_1^2 + d_2^2 + d_3^2 = (d_2 + d_3)^2 - d_1^2 > 0
\end{equation}
and
\begin{equation}
\label{IQ_22}
2d_2d_3+d_1^2-d_2^3-d_3^2 = (d_1-d_2+d_3)(d_1+d_2-d_3) \leq 0.
\end{equation}
From (\ref{IQ_21}) and (\ref{IQ_22}) we can conclude that conjunction
(\ref{KONJUKCIJA_2}) is true and conjunction (\ref{KONJUKCIJA_1}) is not true.
From conjunction (\ref{KONJUKCIJA_2}) follows $d_2^2+d_3^2-d_1^2 \geq 2d_2d_3
> d_1^2 - d_2^2 - d_3^2$ and therefore, $d_2^2 + d_3^2 > d_1^2$, i.e. $\beta_1 > 0$.
The implication (\ref{IMPLIKACIJA_2}) is similarly verified in the case of the inequality
$d_1 + d_2 \leq d_3$.~\stop

\begin{lm}
\label{LEMA_5}
In the metric space $X$ the condition $d_2 + d_3 \leq d_1$ is equivalent to the conjunction
$\alpha_1 \leq 0$ and $\beta_1 \leq 0$.
\end{lm}
{\bf Proof.} $(\Longrightarrow)$ Let for the point $M$ the condition $d_2 + d_3 \leq d_1$
be true. On the basis of equality (\ref{PROIZVOD}) and on the basis of lemma \ref{LEMA_3}
it follows $\alpha_1 \leq 0$. Therefore, on the basis of lemma \ref{LEMA_4}, it follows
$\beta_1 \leq 0$.

\smallskip
\noindent
$(\Longleftarrow)$ Let for the point $M$ conjunction $\alpha_1 \leq 0$ and $\beta_1 \leq 0$
be true. Then from the conjunction
\begin{equation}
\alpha_1 = (d_2 + d_3 - d_1)(d_2 + d_3 + d_1)(2d_2d_3 - \beta_1) \leq 0
\;\;\mbox{and}\;\;
\beta_1 \leq 0
\end{equation}
follows the condition $d_2 + d_3 \leq d_1$.~\stop
\begin{lm}
\label{LEMA_6}
In Ptolemaic metric space $X$ an inequality $\alpha_1^{(1)} \leq 0$
is true iff $b \geq a$ or $c \geq a$.
\end{lm}
{\bf Proof.} On the basis of lemma \ref{LEMA_2} if $a \leq c$ then for the point $B$ we have:
$\alpha_1^{(1)} = a - c \leq 0$ or if $a \leq b$ then for the point $C$ we have:
$\alpha_1^{(1)} = b - a \leq 0$. Conversely, let $a > b,\:\!c$ be true.
Let $M \!\in\! X \backslash \{A,B,C\}$ be any point. Then on the basis of Ptolemaic inequality
\begin{equation}
\label{PT_1}
c \cdot d_3 + b \cdot d_2 \geq a \cdot d_1
\end{equation}
and assumption $a > b,\:\!c$ we can conclude
\begin{equation}
\label{PT_2}
\begin{array}{c}
(c-a)d_3 + (b-a)d_2 + a(d_2+d_3-d_1) \geq 0      \\[2.0 ex]
\Longrightarrow \alpha_1^{(1)} = d_2 + d_3 - d_1 > 0.
\end{array}
\end{equation}
By contraposition the statement follows. ~\stop

\medskip \noindent
On the basis of the previous lemmas we can conclude the following theorem is true.

\begin{th}
\label{TEOREMA_1}
In the metric space $X$ a point $M$ fulfills $\alpha_1^{(1)} = d_2 + d_3 - d_1 \leq 0$
iff  $\alpha_1 \leq 0$ and $\beta_1 \leq 0$ are true. In Ptolemaic metric space
$X$ the set of these  points $M$ is non-empty iff$:$
\begin{equation}
\label{DISJUNKCIJA1}
b \geq a \;\; \mbox{or} \;\; c \geq a.
\end{equation}
\end{th}

\medskip

\leftline{\large \boldmath \bf Inequalities $d_2 + d_3 \leq d_1$,
$d_3 + d_1 \leq d_2$, $d_1 + d_2 \leq d_3$}

\bigskip

\noindent
Let us determine set of  points $M$ in (Ptolemaic) metric spaces for which
some inequalities in (\ref{restart}) are true. With respect to point $A$
we  formed functions (\ref{Alpha}) and (\ref{Beta}). Next, with respect
to point $B$ let us form functions:
\begin{equation}
\alpha_2 = \alpha_2(M) = 4 d_3^2 d_1^2 - {\big (}d_2^2 - (d_3^2 + d_1^2){\big )}^2,
\end{equation}
\begin{equation}
\beta_2 = \beta_2(M) = d_3^2 + d_1^2 - d_2^2
\end{equation}
and with respect to $C$ point let us form functions:
\begin{equation}
\alpha_3 = \alpha_3(M) = 4 d_1^2 d_2^2 - {\big (}d_3^2 - (d_1^3 + d_2^2){\big )}^2,
\end{equation}
\begin{equation}
\beta_3 = \beta_3(M) = d_1^2 + d_2^2 - d_3^2.
\end{equation}
The following equality $\alpha_1 = \alpha_2 = \alpha_3$ is true.
Analogously to the theorem \ref{TEOREMA_1} we can conclude the
following theorems are true.
\begin{th}
\label{TEOREMA_2}
In the metric space $X$ point $M$ fulfills $\alpha_1^{(2)} = d_3 + d_1 - d_2 \leq 0$
iff  $\alpha_1 \leq 0$ and $\beta_2 \leq 0$ are true. In  Ptolemaic metric space
$X$ the set of these  points $M$ is non-empty iff$:$
\begin{equation}
\label{DISJUNKCIJA2}
c \geq b \;\; \mbox{or} \;\; a \geq b.
\end{equation}
\end{th}
\begin{th}
\label{TEOREMA_3}
In the metric space $X$ point $M$ fulfills $\alpha_1^{(3)} = d_1 + d_2 - d_3 \leq 0$
iff $\alpha_1 \leq 0$ and $\beta_3 \leq 0$ are true. In  Ptolemaic metric space
$X$ the set of these  points $M$ is non-empty iff$:$
\begin{equation}
\label{DISJUNKCIJA3}
a \geq c \;\; \mbox{or} \;\; b \geq c.
\end{equation}
\end{th}

\noindent
For (Ptolemaic) metric space $X$ the set of the points $M$ with M\"{o}bius-Pompe\"{i}u
metric property fulfill a conjunction:
\begin{equation}
\label{KONJUKCIJA_ALPHAS}
\alpha_1^{(1)} > 0
\;\;\; \mbox{and} \;\;\;
\alpha_1^{(2)} > 0
\;\;\; \mbox{and} \;\;\;
\alpha_1^{(3)} > 0.
\end{equation}

\noindent
Using theorems \ref{TEOREMA_1}, \ref{TEOREMA_2} and \ref{TEOREMA_3} we can
determine when some inequalities in (\ref{KONJUKCIJA_ALPHAS}) are not true.

\medskip

\noindent
Finally, in the following example let us illustrate a set of points in
$\mbox{\bf R}^2$ with M\"{o}bius-Pompe\"{i}u metric property, with
respect to three fixed points $A, B, C \in \mbox{\bf R}^2$, if we use
metrics $\mbox{\textit{\texttt{d}}}$ and $\overline{\mbox{\textit{\texttt{d}}}
\mathstrut}$ from the example \ref{Example}.

\begin{ex}
{\boldmath $1^{0}.$}
Let in the plane $\mbox{\bf R}^2$ the Euclidean metric $\mbox{\textit{\texttt{d}}}$ is used.
By picture {\rm 1} we illustrate the case of the triangle $ABC$ for which
$\mbox{\textit{\texttt{a}}} > \mbox{\textit{\texttt{c}}} > \mbox{\textit{\texttt{b}}}$
is true. Then $\alpha_1^{(1)}>0$ is true $($the curve $\alpha_1^{(1)}\!=\!0$,
on the basis of the theorem {\rm \ref{TEOREMA_1}}, has empty interior and border$)$,
otherwise the curves $\alpha_1^{(2)}\!=\!0$, $\alpha_1^{(3)}\!=\!0$ have non-empty
interior and border. We can form a non-degenerative triangle from the remaining points.


\vspace*{120.0 mm}

\centerline{\mbox{\rm Picture 1.}}

\newpage

\noindent
In the case of the equilateral triangle $ABC$ the curves $\alpha_{1}^{(1)}\!=\!0$,
$\alpha_{1}^{(2)}\!=\!0$ and $\alpha_{1}^{(3)}\!=\!0$ transform onto the $($smaller$)$
arcs $\arc{B}{C}$, $\arc{C}{A}$ and $\arc{A}{B}$ of the circumcircle. Hence, we have
M\"{o}bius-Pompe\"{i}u theorem in the following form$:$ for equilateral triangle $ABC$
the set of  points $M$ in the plane, such that from distances $\mbox{\textit{\texttt{d}}}_1
= \mbox{\textit{\texttt{d}}}(M,A)$, $\mbox{\textit{\texttt{d}}}_2
= \mbox{\textit{\texttt{d}}}(M,B)$ and $\mbox{\textit{\texttt{d}}}_3
= \mbox{\textit{\texttt{d}}}(M,C)$ one can form a degenerative triangle,
is circumcircle; from the other points in the plane we can form non-degenerative
triangle.

\medskip
\noindent
{\boldmath $2^{0}.$}
Let in the plane $\mbox{\bf R}^2$ the chordal metric $\overline{\mbox{\textit{\texttt{d}}}\mathstrut}$
is used. Let $A,B,C \!\in\! \mbox{${\cal S} \backslash \{(0,0,1)\}$}$ be points
on the
unit
Riemann sphere ${\cal S}$, with uniquely determined projections$:$
$$
A^{'} = {\cal P}^{-1}(A) = a_1 \!\,+\!\, b_1 i,\,
B^{'} = {\cal P}^{-1}(B) = a_2 \!\,+\!\, b_2 i,\,
C^{'} = {\cal P}^{-1}(C) = a_3 \!\,+\!\, b_3 i \in \mbox{\bf C}
$$
with inversely stereographical projection
from the north pole$:$
$$
{\cal P}^{-1}
=
{\cal P}^{-1}(x,y,z)
=
{\Big (} \displaystyle\frac{x}{1-z} {\Big )}
+
{\Big (} \displaystyle\frac{y}{1-z} {\Big )} i
:
{\cal S} \backslash \{(0,0,1)\} \longrightarrow \mbox{\bf C}.
$$

\noindent
\begin{minipage}[b]{80.0 mm}
Through points $A,B,C$ on the Riemann sphere let us set great circles $($picture {\rm 2}$)$.
In the complex plane we uniquely determine images of great circles as corresponding  circles
through points $A^{'},B^{'},C^{'}$ \mbox{$($picture~{\rm 3}$)$}. By picture~{\rm 3}
we illustrate the case of points $A^{'},B^{'},C^{'}$ for which
$\overline{\mbox{\textit{\texttt{b}}}\mathstrut}
\!\,>\!\,
\overline{\mbox{\textit{\texttt{c}}}\mathstrut}
\!\,>\!\,
\overline{\mbox{\textit{\texttt{a}}}\mathstrut}$
and $\overline{\mbox{\tt k}\mathstrut} \neq 0$ are
true.
Then $\overline{\alpha\mathstrut}_1^{(2)}\!>\!0$ $($the curve
$\overline{\alpha\mathstrut}_1^{(2)}\!=\!0$,
on the basis of the theorem {\rm \ref{TEOREMA_2}}, has empty interior and border$)$,
otherwise curves \mbox{$\overline{\alpha\mathstrut}_1^{(1)}\!=\!0$},
\mbox{$\overline{\alpha\mathstrut}_1^{(3)}\!=\!0$} have non-empty
interior and border. From the remaining points we can form a non-degenerative triangle.
\end{minipage}
\hspace*{19.0 mm}{\mbox{\rm Picture 2.}}


\vspace*{65.0 mm}


\centerline{\mbox{\rm Picture 3.}}

\break

\noindent
Let us consider the case when $A$, $B$, $C$ are chordally equidistantly arranged points
on the Riemann sphere ${\cal S}$. Then the set of  points $M$ on the Riemann sphere,
being such that from chordal distances
$\overline{\mbox{\textit{\texttt{d}}}\mathstrut}_1
= \overline{\mbox{\textit{\texttt{d}}}\mathstrut}(M,A)$,
$\overline{\mbox{\textit{\texttt{d}}}\mathstrut}_2
= \overline{\mbox{\textit{\texttt{d}}}\mathstrut}(M,B)$
and
$\overline{\mbox{\textit{\texttt{d}}}\mathstrut}_3
= \overline{\mbox{\textit{\texttt{d}}}\mathstrut}(M,C)$
one can form a degenerative triangle, is circumcircle; from other points on the Riemann
sphere one can form a non-degenerative triangle.
Using inverse stereographical projection ${\cal P}^{-1}$
we can conclude that analogous statement in complex plane
$\mbox{\bf C}$ is valid if we use chordal metric
$\overline{\mbox{\textit{\texttt{d}}}\mathstrut}$.
\end{ex}

\newpage


\end{document}